\RequirePackage{fix-cm}
\documentclass[smallcondensed]{svjour3}     
\smartqed  
\usepackage{graphicx}
\usepackage{amssymb}
\usepackage{algorithm}
\usepackage{caption}
\usepackage{latexsym}
\usepackage{mathtools}
\usepackage{booktabs}
\usepackage{enumerate}
\usepackage{tikz}
\usepackage{lmodern}
\usepackage{microtype}
\usepackage{subfigure}
\DeclarePairedDelimiter{\norm}{ \lVert }{ \rVert }

\DeclareMathOperator{\supp}{supp}

\begin{document}

\title{Partial Gradient Optimal Thresholding Algorithms for a Class of Sparse Optimization Problems
\thanks{This work was founded by the National Natural Science Foundation of China (NSFC) under the grant 12071307.} }
\titlerunning{Partial Gradient Optimal $k$-thresholding}

\author{Nan Meng \and Yun-Bin Zhao \and Michal Ko\v{c}vara}


\institute{Nan Meng \at
              School of Mathematics, University of Birmingham, Edgbaston, Birmingham B15 2TT, United Kingdom \\
              \email{nxm563@bham.ac.uk}
           \and
             Yun-Bin Zhao \at
              Shenzhen Research Institute of Big Data, Chinese University of Hong Kong, Shenzhen, China \\
              \email{yunbinzhao@cuhk.edu.cn}
             \and
             Michal Ko\v{c}vara \at
             School of Mathematics, University of Birmingham, Edgbaston, Birmingham B15 2TT, United Kingdom \\
             \email{m.kocvara@bham.ac.uk}
}

\date{Received: date / Accepted: date}

\maketitle

\begin{abstract}
The optimization problems with a sparsity constraint is a class of important global optimization problems.
A typical type of thresholding algorithms for solving such a problem adopts the traditional full steepest descent direction or Newton-like direction as a search direction to generate an iterate on which a certain thresholding is performed.
Traditional hard thresholding discards a large part of a vector when the vector is dense.
Thus a large part of important information contained in a dense vector has been lost in such a thresholding process.
Recent study [Zhao, SIAM J Optim, 30(1), pp. 31-55, 2020] shows that the hard thresholding should be applied to a compressible vector instead of a dense vector to avoid a big loss of information.
On the other hand, the optimal $k$-thresholding as a novel thresholding technique may overcome the intrinsic drawback of hard thresholding, and performs thresholding and objective function minimization simultaneously.
This motivates us to propose the so-called partial gradient optimal thresholding method in this paper, which is an integration of the partial gradient and the optimal $k$-thresholding technique.
The solution error bound and convergence for the proposed algorithms have been established in this paper under suitable conditions.
Application of our results to the sparse optimization problems arising from signal recovery is also discussed.
Experiment results from synthetic data indicate that the proposed algorithm called PGROTP is efficient and comparable to several existing algorithms.
\keywords{Sparse optimization \and Signal recovery \and Optimal $k$-thresholding \and Partial gradient method \and Error bound \and Convergence}
\end{abstract}

\section{Introduction}
The optimization problem with a sparsity constraint
\begin{equation} \label{prob}
    \min_{x \in \mathbb{R}^n} \left\{f(x) : ~ \norm{x}_0 \le k \right\}
\end{equation}
arises in many practical fields such as compressive sensing \cite{eldar2012compressed,foucart2013mathematical,zhao2018sparse}, signal processing \cite{boche2019compressed}, wireless communication \cite{choi2017compressed}, pattern recognition \cite{patel2011sparse}, to name a few.
In the model \eqref{prob}, $\|\cdot\|_0$ is often called the $\ell_0$-norm which counts the number of nonzero entries of a vector.
Depending on application, the function $f(x)$ may take different specific forms.
For instance, in compressed sensing scenarios, $f(x)$ is usually taken as $\left\|y-Ax\right\|_2^2$ which is an error metric for signal measurements.
The problem \eqref{prob} is known to be NP-hard, and the main difficulty for solving this problem lies in locating the position of nonzero entries of a feasible sparse vector at which $f(x)$ is minimized.

The algorithms for solving \eqref{prob} can be sorted into several categories including convex optimization methods, heuristic algorithms, and thresholding algorithms.
The convex optimization methods include $\ell_1$-minimization \cite{chen2001atomic,candes2005decoding}, reweighted $\ell_1$-minimization \cite{candes2008enhancing,zhao2012reweighted}, and dual-density-based reweighted $\ell_1$-minimization \cite{zhao2015new,zhao2017constructing,zhao2018sparse}.
The heuristic-type methods include orthogonal matching pursuit (OMP) \cite{tropp2007signal,needell2010signal,cai2011orthognal}, compressive sampling matching pursuit (CoSaMP) \cite{needell2009cosamp}, subspace pursuit (SP) \cite{dai2008subspace,dai2009subspace}, and their variants.
Thresholding-type algorithms attract much attention due to their easy implementation and low computational complexity \cite{blumensath2009iterative,blumensath2010normalized,foucart2011hard,blumensath2012accelerated,bouchot2016hard,khanna2018iht,zhao2020optimal}.

The key step for thresholding-type iterative algorithms can be stated as
\begin{equation} \label{threshold}
	x^{p+1} = {\cal T}_k \left(x^p+\lambda d\right),
\end{equation}
where ${\cal T}_k (\cdot)$ represents a thresholding operator that is used to produce a $k$-sparse vector, $\lambda$ denotes the stepsize and $d$ is a search direction at the current iterate $x^p$.
Throughout the paper, a vector $x$ is said to be $k$-sparse if $\norm{x}_0 \le k$.
Several thresholding operators are widely used in the literature, such as the hard thresholding \cite{blumensath2009iterative,blumensath2010normalized,foucart2011hard,blumensath2012accelerated,bouchot2016hard,khanna2018iht}, soft thresholding \cite{donoho1995denoising,fornasier2008iterative,huang2013soft,liu2016projected}, and optimal $k$-thresholding \cite{zhao2020optimal,zhao2020analysis}.
The steepest descent direction \cite{blumensath2008iterative,blumensath2009iterative,garg2009gradient,blumensath2010normalized,bouchot2016hard,nan2020newton} and Newton-type direction \cite{yuan2014newton,nan2021newton,zhou2021global,zhou2020subspace} are two search directions that are used by many researchers.

Let ${\cal H}_k$ denote the hard thresholding operator which retains the largest $k$ magnitudes and zeroes out other entries of a vector.
By setting ${\cal T}_k = {\cal H}_k$ and $d = -\nabla f(x)$, the iterative formula \eqref{threshold} is reduced to
\begin{equation} \label{IHTite}
	x^{p+1} = {\cal H}_k \left(x^p - \lambda \nabla f(x^p) \right),
\end{equation}
where $\nabla f (x^p)$ is the gradient of $f$ at $x^p$.
The formula \eqref{IHTite} is the well-known iterative hard thresholding (IHT) initially studied in \cite{blumensath2008iterative,blumensath2009iterative}.
The IHT can be enhanced by either attaching an orthogonal projection (a pursuit step) to obtain the so-called hard thresholding pursuit (HTP) method \cite{foucart2011hard,bouchot2016hard} or by using an adaptive stepsize strategy to yield the so-called normalized iterative hard thresholding (NIHT) \cite{blumensath2010normalized}.
While the algorithm \eqref{IHTite} can reconstruct the vector under suitable conditions (see \cite{blumensath2009iterative,foucart2013mathematical,zhao2020improved}), but as pointed in \cite{zhao2020optimal,zhao2020analysis}, the operator ${\cal H}_k$ may cause certain numerical problems as well.

To improve the performance of IHT, Zhao \cite{zhao2020optimal,zhao2020analysis} recently proposed the optimal $k$-thresholding technique which stresses that thresholding of a vector should be performed simultaneously with objective function reduction in the course of iterations.
Replacing ${\cal H}_k$ by the optimal $k$-thresholding operator ${\cal Z}^{\#}_k$ in \eqref{IHTite} leads to the following iterative optimal $k$-thresholding scheme:
\begin{equation*}
	x^{p+1} = {\cal Z}^{\#}_k \left(x^p - \lambda \nabla f(x^p) \right).
\end{equation*}
The optimal $k$-thresholding of a vector $u \in \mathbb{R}^n$ with respect to the objective function $f(x)$ is defined as ${\cal Z}^{\#}_k (u) := u \otimes w^*$ with
\begin{equation} \label{OTbinary}
	w^{*}=\arg \min _{w}\left\{f(u \otimes w): ~ \mathbf{e}^{T} w=k, ~ w \in\{0,1\}^{n}\right\} ,
\end{equation}
where $\mathbf{e} \in \mathbb{R}^{n}$ is the vector of ones and $\otimes$ denotes the Hadamard product of two vectors.
To avoid solving the above binary optimization problem, Zhao \cite{zhao2020optimal} suggests solving the tightest convex relaxation of \eqref{OTbinary} instead.
That is, replacing the binary constraint by its convex hull, we obtain the following convex relaxation of the problem \eqref{OTbinary}:
\begin{equation*}
	\overline{w}=\arg \min _{w}\left\{f(u \otimes w): ~ \mathbf{e}^{T} w=k, ~ 0 \le w \le \mathbf{e} \right\}.
\end{equation*}
The vector $u \otimes \overline{w}$ is called the relaxed optimal $k$-thresholding of $u$.

The hard thresholding operator in \eqref{IHTite} discards a large important part of a vector when the vector is dense.
This means some important information of the vector has been lost in the process of hard thresholding.
As pointed out in \cite{zhao2020optimal,zhao2020analysis}, the hard thresholding should be applied to a compressible vector instead of a dense vector in order to avoid a big loss of information.
Note that the vector $x^p - \lambda \nabla f(x^p)$ in \eqref{IHTite} is usually dense since the search direction $-\nabla f(x^p)$ is not necessarily sparse.
This motivates us to adopt the partial gradient instead of the full gradient as a search direction in order to generate the following sparse or compressible vector:
\[ u^{p} := x^p - \lambda {\cal H}_q (\nabla f(x^p)) \]
on which some thresholding is then performed to generate the next iterate $x^{p+1}$.
In the above formula, the integer number $q>0$ controls the number of elements extracted from the full gradient.
In other words, we only use $q$ significant entries of the gradient as our search direction.
We may use the hard thresholding of $u^p$ to produce an iterate satisfying the constraint of the problem \eqref{prob}.
However, as we pointed out before, the optimal $k$-thresholding is more powerful and more efficient than the hard thresholding.
This stimulates the following iterative scheme:
\begin{equation} \label{PGOTkey}
	x^{p+1} = {\cal Z}^{\#}_k \left(x^p - \lambda {\cal H}_q (\nabla f(x^p))\right).
\end{equation}
This is refer to as the partial gradient optimal thresholding (PGOT) algorithm in this paper, which is described in detail in Section 2.
The enhanced version of PGOT is called the partial gradient relaxed optimal thresholding (PGROT).
In order to maintain the $k$-sparsity of the iterate, a further enhancement of PGOT and PGROT can be made by adding a pursuit step to PGROT to eventually obtain the more efficient algorithm called the partial gradient relaxed optimal thresholding pursuit (PGROTP), which is treated as the final version of the proposed algorithm actually being used to solve the problems.
The solution error bound and convergence analysis for our algorithms with $q$ in the range $q \in [2k,n]$ are shown under the assumption of restricted isometry property (RIP).
Simulations from synthetic data indicate that PGROTP algorithm is robust and comparable to several existing methods.

The paper is organized as follows.
The algorithms are described in Section 2.
The error bounds and global convergence of the proposed algorithms are established in Section 3.
Numerical results are given in Section 4 and conclusions are given in Section 5.

\subsection{Notations}
We first introduce some notations used throughout the paper.
$\mathbb{R}^n$ is the $n$-dimensional Euclidean space, and $\mathbb{R}^{m \times n}$ is the set of $m \times n$ matrices.
Vector $\mathbf{e}$ is the vector of ones.
Denote $[N]$ as the set $\{1,\dots,n\}$.
Given a set $\Omega \subseteq [N]$, $\overline{\Omega} \coloneqq [N] \backslash \Omega$ denotes the complement set of $\Omega$, and $|\Omega|$ is the cardinality of set $\Omega$.
For a vector $x \in \mathbb{R}^n$, $x_{\Omega} \in \mathbb{R}^n$ denotes the vector obtained from $x$ by retaining elements indexed by $\Omega$ and zeroing out the remaining ones.
The set $\supp(x) = \{x_i, i \neq 0\}$ is called the support of $x$, ${\cal L}_k (x)$ denotes the support of ${\cal H}_k (x)$, and ${\cal Z}^{\#}_k (\cdot)$ denotes the optimal $k$-thresholding operator.
For a matrix $A$, $A^T$ denotes its transpose.
The notation $\otimes$ represents the Hadamard product of two vectors, i.e., $u \otimes v = [u_1v_1, \dots, u_nv_n]^T$.
Given a number $\alpha$, $\lceil \alpha \rceil$ is the smallest integer number that is larger than or equal to $\alpha$.

\section{Algorithms}
In this paper, we focus on the following specific objective function:
\begin{equation} \label{function}
	f(x) = \frac{1}{2} \norm*{y-Ax}_2^2,
\end{equation}
where $A$ is a given $m \times n$ matrix with $m \ll n$, and $y$ is a given vector. Using this quadratic function, the model \eqref{prob} becomes
\begin{equation} \label{problem}
	\min_{x \in \mathbb{R}^n} \left\{ \frac{1}{2} \norm*{y-Ax}_2^2 : \norm{x}_0 \le k \right\}.
\end{equation}
This problem has been widely used in signal reconstruction via compressive sensing and in many other application settings.
The gradient of the function \eqref{function} is given as
$$\nabla f(x) = -A^T(y-Ax).$$
By using the major part of this specific gradient, we define the $u^p$ as follows:
\begin{align*}
	u^p = x^p + \lambda {\cal H}_{q} \left( A^T(y-Ax^p) \right),
\end{align*}
where $0 < q \le n$ is an integer number.
The optimal thresholding  of $u^p$ with respect to the function \eqref{function} is given by ${\cal Z}^{\#}_k (u^p) = u^p \otimes w^*$, where
\begin{equation*}
	w^{*}=\arg \min _{w}\left\{\norm*{y-A (u \otimes w)}_2^2: ~ \mathbf{e}^{T} w=k, ~ w \in\{0,1\}^{n}\right\} .
\end{equation*}
Thus the partial gradient optimal $k$-thresholding (PGOT) algorithm \eqref{PGOTkey} for solving problem \eqref{problem} can be stated as
\[ x^{p+1} = {\cal Z}_k^{\#} \left( x^p + \lambda {\cal H}_q \left( A^T(y-Ax^p)\right) \right). \]
For simplicity of algorithmic description and analysis, we set $\lambda = 1$ throughout the rest of the paper.
It should be stressed that in practical applications, suitable stepsize should be used in order to speed up the convergence of the algorithms.
By the definition of ${\cal Z}_k^{\#}$, the PGOT algorithm can be described explicitly as follows.

\begin{algorithm}
\caption{Partial Gradient Optimal $k$-Thresholding (PGOT)} \label{PGOT}
\begin{itemize}
    \item Input: matrix $A$, vector $y$, sparsity level $k$, integer number $q \ge k$, and initial point $x^0 = 0$.
    \item Iteration:
    \begin{align}
        u^p & = x^p + {\cal H}_{q} \left( A^T(y-Ax^p) \right) , \nonumber \\
        w^p & = \arg \min_w \{ \norm*{y - A (w \otimes u^p)}_2^2 : ~ \textbf{e}^T w = k, ~ w \in \{0,1\}^n \} , \tag{OP} \\
        x^{p+1} & = w^p \otimes u^p ~~ (={\cal Z}_k^{\#} (u^p)). \nonumber
    \end{align}
    \item Output: $k$-sparse vector $\hat{x}$.
\end{itemize}
\end{algorithm}

To avoid solving the integer programming problem (OP), as suggested in \cite{zhao2020optimal}, the binary constraint in (OP) $w \in \{0,1\}^n$ can be relaxed to $0 \le w \le \mathbf{e}$ so that we obtain the partial gradient relaxed optimal thresholding (PGROT) algorithm.

\begin{algorithm}[H]
\caption{Relaxed Partial Gradient Optimal $k$-Thresholding (PGROT)} \label{PGROT}
\begin{itemize}
    \item Input: matrix $A$, vector $y$, sparsity level $k$, integer number $q \ge k$, and initial point $x^0 = 0$.
    \item Iteration:
    \begin{align*}
        u^p & = x^p + {\cal H}_{q} \left( A^T(y-Ax^p) \right), \\
        \overline{w}^p & = \arg \min_w \{ \norm*{y - A (w \otimes u^p)}_2^2 : ~ \textbf{e}^T w = k, ~ 0 \le w \le \textbf{e} \}, \tag{ROT} \\
        x^{p+1} & = {\cal H}_k (\overline{w}^p \otimes u^p) .
    \end{align*}
    \item Output: $k$-sparse vector $\hat{x}$.
\end{itemize}
\end{algorithm}

The solution $\overline{w}^p$ to (ROT) is not $k$-sparse in general.
So the purpose of the final thresholding step in PGROT to restore the $k$-sparsity of iterate.
It is worth emphasizing that the use of ${\cal H}_k$ here is different from the settings in traditional IHT, since the vector $u^p$ generated by the partial gradient is $(k+q)$-sparse instead of being a usually dense vector in IHT.

The PGROT can be further enhanced by including a pursuit step (i.e., an orthogonal projection step) to find a possibly better iterate than the point generated by PGROT.
This consideration leads to so-called PGROTP algorithm which is described in Algorithm 3.
In the next section, we perform a theoretical analysis for the proposed algorithms focusing on establishing the error bound for the solution of the problem.

\begin{algorithm}
\caption{Relaxed Partial Gradient Optimal $k$-thresholding Pursuit (PGROTP)}
\begin{itemize}
    \item Input: matrix $A$, vector $y$, sparsity level $k$, integer number $q \ge k$, and initial point $x^0 = 0$.
    \item Iteration:
    \begin{align*}
        u^p & = x^p + {\cal H}_{q} \left( A^T(y-Ax^p) \right), \\
        \overline{w}^p & = \arg \min_w \{ \norm*{y - A (w \otimes u^p)}_2^2 : ~ \textbf{e}^T w = k, ~ 0 \le w \le \textbf{e} \}, \\
        S^{p+1} & = {\cal L}_k (\overline{w}^p \otimes u^p), \\
        x^{p+1} & = \arg \min_z \{ \norm*{y-Az}_2 : ~ \supp(z) \subseteq S^{p+1} \} . \tag{Pursuit step}
    \end{align*}
    \item Output: $k$-sparse vector $\hat{x}$.
\end{itemize}
\end{algorithm}

\section{Error Bound and Convergence Analysis}

In this section, we establish the error bounds for the solution of the problem via the proposed algorithms.
The purpose is to estimate the distance between the iterate $x^p$, generated by the proposed algorithms, and the global solution of the problem \eqref{problem}.
As an implication of the error bounds, the global convergence of our algorithms can be instantly obtained for the problem \eqref{problem} arising from the scenarios of sparse signal recovery.

Before going ahead, we first introduce the restricted isometry constant (RIC) of the matrix $A$.

\begin{definition}\cite{candes2005decoding,foucart2013mathematical}
The $s$-th order restricted isometry constant (RIC) $\delta_s$ of a matrix $A \in \mathbb{R}^{m \times n}$ is the smallest number $\delta_s \ge 0$ such that
\[  \left(1-\delta_{s}\right)\|x\|_{2}^{2} \leq\|A x\|_{2}^{2} \leq\left(1+\delta_{s}\right)\|x\|_{2}^{2}  \]
for all $s$-sparse vector $x$, where $s>0$ is an integer number.
\end{definition}

It is usually to say that matrix $A$ satisfies the RIP of order $s$ if $\delta_s < 1$.
It is well known that the random matrices such as Bernoulli and Gaussian matrices are widely used in applications as they satisfy the RIP with an overwhelming probability \cite{candes2005decoding,candes2006robust}.
The following property of RIC is frequently used in our paper.
\begin{lemma} \emph{\cite{candes2005decoding,needell2009cosamp,foucart2011hard}} \label{RIP}
Suppose matrix $A$ satisfy the RIP of order $k$.
Given a vector $u \in \mathbb{R}^n$ and a set $\Omega \subseteq [N]$, one has
\begin{enumerate}[(i)]
\item  $\norm*{\left(\left(I-A^{T} A\right) v\right)_{\Omega}}_{2} \leq \delta_{t}\|u\|_{2} $ if $|\Omega \cup \supp(v)| \leq t$.
\item $\norm*{(A^T u)_{\Omega}}_2 \le \sqrt{1+\delta_t} \norm*{u}_2$ if $|\Omega| \le t$.
\end{enumerate}
\end{lemma}

\subsection{Main results for PGOT} \label{section3}
The following two technical results concerning the properties of optimal $k$-thresholding and hard thresholding operators are useful.

\begin{lemma} \emph{\cite{zhao2020optimal,zhao2020analysis}} \label{lem2}
Let $x^* \in \mathbb{R}^n$ be the solution to the problem \eqref{problem} and denote by $\eta = y-Ax^*$.
Given an arbitrary vector $u \in \mathbb{R}^n$, let ${\cal Z}^{\#}_k (u)$ be the optimal $k$-thresholding vector of $u$.
Then for any $k$-sparse binary vector $w^* \in \{0,1\}^n$ satisfying $\operatorname{supp}(x^*) \subseteq \operatorname{supp}(w^*)$, one has
\begin{equation*}
	\norm*{{\cal Z}^{\#}_k (u) - x^*}_2 \le \sqrt{\frac{1+\delta_k}{1-\delta_{2k}}} \norm*{(x^* - u) \otimes w^*}_2 + \frac{2}{\sqrt{1-\delta_{2k}}} \| \eta \|_2.
\end{equation*}
\end{lemma}
This result can be found from the proof of Theorem 4.3 in \cite{zhao2020optimal}.
\begin{lemma} \emph{\cite{zhao2020improved}} \label{hard}
Let $q \ge s$ be an integer number.
For any vector $z \in \mathbb{R}^{n}$ and any $s$-sparse vector $u \in \mathbb{R}^{n}$, one has
    \[
    \left\|u-\mathcal{H}_{q}(z)\right\|_{2} \le \frac{\sqrt{5}+1}{2} \|(u-z)_{\Lambda \cup \Omega}\|_{2},
    \]
    where $\Lambda=\supp(u)$ and $\Omega =\supp\left(\mathcal{H}_{q}(z)\right)$.
\end{lemma}
When $s \le q$, a $s$-sparse vector is also $q$-sparse.
Thus Lemma \ref{hard} above follows exactly from Lemma 2.2 in \cite{zhao2020improved}.
We are ready to prove the error bound and global convergence of the algorithm PGOT.
\begin{theorem}
Let $x^* \in \mathbb{R}^n$ be the solution of the problem \eqref{problem} and $\eta := y-Ax^*$.
Let $q \geq 2 k$ be a positive integer number.
Suppose the restricted isometry constant of $A$ satisfies
\begin{equation} \label{thm1-1}
	\delta_{3k} < \alpha^*,
\end{equation}
where $\alpha^* \in (0,1)$ is the unique real root of the univariate equation
\begin{equation} \label{alpha}
	\left(\left\lceil \frac{q}{k} \right\rceil+1\right)^{2} \alpha^{3}+\left(\left\lceil \frac{q}{k} \right\rceil+1\right)^{2} \alpha^{2}+\frac{2}{3+\sqrt{5}} \alpha -\frac{2}{3+\sqrt{5}} = 0.
\end{equation}
Then the sequence $\left\{x^{p}\right\}$ generated by PGOT satisfies that
\begin{equation*}
\norm*{x^{p+1}-x^*}_2 \le \rho^p \norm*{x^{0}-x^*}_2 + \frac{\tau}{1-\rho} \norm{\eta}_2,
\end{equation*}
where
\begin{equation} \label{1-8}
\rho := \frac{\sqrt{5}+1}{2} \left(\left\lceil\frac{q}{k}\right\rceil+1\right) \delta_{3 k} \sqrt{\frac{1+\delta_{k}}{1-\delta_{2 k}}} <1
\end{equation}
is guaranteed under the condition \eqref{thm1-1}, and the constant $\tau$ is given as
\begin{equation} \label{1-9}
\tau = \frac{(\sqrt{5}+1) \left(\left\lceil\frac{q}{k}\right\rceil+1\right) (1+\delta_{2k})+4}{2\sqrt{1-\delta_{2k}}} .
\end{equation}
In particular, if $\eta=0$, the sequence $\left\{x^{p}\right\}$ generated by PGOT converges to $x^*$.
\end{theorem}
\emph{Proof.}
Let $\eta = y-Ax^*$ and $u^p, x^{p+1}$ be the vectors generated at $p$-th iteration of PGOT, i.e., $u^{p}=x^{p}+\mathcal{H}_{q}\left(A^{T}\left(y-A x^{p}\right)\right)$ and $x^{p+1} = {\cal Z}^{\#}_k (u^p)$.
Let $w^* \in \{0,1\}^n$ be a $k$-sparse vector such that $\supp(x^*) \subseteq \supp(w^*)$.
Applying Lemma \ref{lem2} leads to
\begin{align}
\norm*{x^*-x^{p+1}}_2 & = \norm*{x^* - {\cal Z}^{\#}_k (u^p)}_2 \nonumber \\
& \leq \sqrt{\frac{1+\delta_{k}}{1-\delta_{2 k}}}\left\|\left(x^*-u^{p}\right) \otimes w^* \right\|_{2}+\frac{2}{\sqrt{1-\delta_{2 k}}}\left\|\eta\right\|_{2} \nonumber \\
& \le \sqrt{\frac{1+\delta_{k}}{1-\delta_{2 k}}}\left\|x^*-u^{p} \right\|_{2}+\frac{2}{\sqrt{1-\delta_{2 k}}}\left\|\eta\right\|_{2} ~~~ \text{(since $0 \le w^* \le \textbf{e}$).} \label{1-7}
\end{align}
Denote $\Omega := {\cal L}_q (A^{T}(y-A x^{p}) )$.
It is easy to see that $|\supp(x^*-x^p)| \le 2k$.
Thus, if $q \ge 2k$, by Lemma \ref{hard}, we have
\begin{align}
\left\|x^*-u^{p}\right\|_{2}
& = \left\|x^*-x^{p}-\mathcal{H}_{q}\left(A^{T}\left(y-A x^{p}\right)\right)\right\|_{2} \nonumber \\
& \le \frac{\sqrt{5}+1}{2} \left\|\left(x^*-x^{p}-A^{T}\left(y-A x^{p}\right) \right)_{\Omega \cup (S \cup S^p)}\right\|_{2} , \label{1-4}
\end{align}
where $S = \supp(x^*)$ and $S^p = \supp(x^p)$.
Given a vector $v \in \mathbb{R}^n$ and two support sets $\Lambda_1, \Lambda_2 \in [N]$.
It is easy to verify $\norm*{v_{\Lambda_1 \cup \Lambda_2}}_2 \le \norm*{v_{\Lambda_1}}_2 + \norm*{v_{\Lambda_2}}_2$.
Therefore,
\begin{align}
\|(x^*-x^{p}-A^{T} &(y-A x^{p}) )_{\Omega \cup (S \cup S^p)}\|_{2} \nonumber
\le \left\|\left(x^*-x^{p}-A^{T}\left(y-A x^{p}\right) \right)_{\Omega}\right\|_{2} \\
& + \left\|\left(x^*-x^{p}-A^{T}\left(y-A x^{p}\right) \right)_{S \cup S^p} \right\|_{2} . \label{1-3}
\end{align}
Noting that $|\supp(x_S-x^p)| \le |S \cup S^p| \le 2k$.
The second term on the right-hand side of \eqref{1-3} can be bounded.
In fact, by Lemma \ref{RIP}, we have
\begin{align}
\|(x^*-x^{p}&-A^{T}(y-A x^{p}) )_{S \cup S^p}\|_{2}  \nonumber \\
& = \norm*{\left((I-A^TA) (x^*-x^p) + A^T\eta\right)_{S \cup S^p}}_2 \nonumber \\
& \le \norm*{\left((I-A^TA) (x^*-x^p)\right)_{S \cup S^p}}_2 + \norm*{(A^T\eta)_{S \cup S^p}}_2 \nonumber \\
& \le \delta_{2k} \left\|x^*-x^{p}\right\|_{2} + \sqrt{1+\delta_{2k}} \left\|\eta\right\|_{2} . \label{1-5}
\end{align}
Setting $t = \left\lceil \frac{q}{k} \right\rceil$, the set $\Omega$ can be separated into $t$ disjoint sets such that $\Omega = T_1 \cup T_2 \dots, T_t$, where $|T_i| \le k$ for $i=1,\dots,t$, and $T_i \cap T_j = \emptyset$ if $i \neq j$.
Thus we have
\begin{align}
\|(x^*-x^{p}&-A^{T}(y-A x^{p}))_{\Omega}\|_{2} \nonumber \\
& \le \sum_{i=1}^{t}\left\|\left(x^*-x^{p}-A^{T}\left(y-A x^{p}\right) \right)_{T_i}\right\|_{2} \nonumber \\
& \le \sum_{i=1}^{t} \left\| \left[ (I-A^TA)(x^*-x^p) \right]_{T_i} \right\|_{2} + \sum_{i=1}^{t} \left\| \left[ A^T\eta\right]_{T_i}\right\|_{2} \nonumber \\
& \le t \delta_{3k} \left\|x^*-x^p\right\|_{2} + t \sqrt{1+\delta_k} \left\|\eta\right\|_2 , \label{1-6}
\end{align}
where the last inequality follows from Lemma \ref{RIP} because of $|T_i \cup \supp(x^*-x^p)| \le 3k$.
Since $\delta_{2k} \le \delta_{3k}$, combining \eqref{1-4}-\eqref{1-6} leads to
\begin{align}
\left\|x^*-u^{p}\right\|_{2}
& \le \frac{\sqrt{5}+1}{2} (t+1) \delta_{3k} \left\|x^*-x^{p}\right\|_{2} + \frac{\sqrt{5}+1}{2}(t+1) \sqrt{1+\delta_{2k}}\left\|\eta\right\|_{2} . \label{1-2}
\end{align}
Substituting \eqref{1-2} into \eqref{1-7} yields
\begin{equation} \label{1-13}
	\left\|x^{p+1}-x^*\right\|_{2} \leq \rho\left\|x^{p}-x^*\right\|_{2}+\tau\left\|\eta\right\|_{2},
\end{equation}
where $\rho$ and $\tau$ are given as \eqref{1-8} and \eqref{1-9}, respectively.
Since $\delta_k \le \delta_{2k} \le \delta_{3k}$, the constant $\rho<1$ is ensured if
\begin{equation} \label{1-10}
\frac{\sqrt{5}+1}{2}(t+1) \delta_{3 k} \sqrt{\frac{1+\delta_{3k}}{1-\delta_{3k}}} < 1.
\end{equation}
Squaring both sides of \eqref{1-10} and rearranging terms yield
\[ g(\delta_{3k}) := (t+1)^2 \delta_{3k}^3 + (t+1)^2 \delta_{3k}^2 + \frac{2}{3+\sqrt{5}} \delta_{3k} - \frac{2}{3+\sqrt{5}} < 0. \]
The gradient of $g$ with respect to $\delta_{3k}$ is given as
\[ \nabla g(\delta_{3k}) = 3(t+1)^2 \delta_{3k}^2 + 2(t+1)^2\delta_{3k} + \frac{2}{3+\sqrt{5}}> 0. \]
Thus the function $g$ is strictly and monotonically increasing over the interval $\delta_{3k} \in (0,1]$.
Note that
\[ g(0) = - \frac{2}{3+\sqrt{5}} < 0 ~~\text{and}~~ g(1) = 2(t+1)^2 > 0. \]
Thus there exists a unique real root $\alpha^*$ for the equation $g(\alpha^*)=0$ in $[0,1]$.
Therefore, $\delta_{3k} < \alpha^*$ ensures that the constant $\rho<1$ in \eqref{1-13}, and hence it follows from \eqref{1-13} that
\begin{equation*}
	\left\|x^{p+1}-x^{*}\right\|_{2} \leq \rho^{p}\left\|x^{0}-x^{*}\right\|_{2}+\frac{\tau}{1-\rho}\|\eta\|_{2},
\end{equation*}
which is exactly the desired error bound.
In particular, when $\eta=0$, it follows immediately from the above error bound that the sequence $\{x^p\}$ generated by PGOT converges to $x^*$ as $p \to \infty$.
\hfill $\Box$

A more explicitly given RIC bound than \eqref{thm1-1} for PGOT can be derived as follows.
Since $\sqrt{\frac{1+\delta_{3k}}{1-\delta_{3k}}} < \frac{1+\delta_{3k}}{1-\delta_{3k}}$, the inequality \eqref{1-10} is guaranteed provided the following inequality is satisfied:
\[ \frac{\sqrt{5}+1}{2} \left(\left\lceil\frac{q}{k}\right\rceil+1 \right) \delta_{3 k} \frac{1+\delta_{3k}}{1-\delta_{3k}} < 1, \]
which can be written as
\begin{equation} \label{1-12}
	\phi \delta_{3k}^2 + \left( \phi + 1 \right)\delta_{3k} - 1 < 0,
\end{equation}
where
\[ \phi = \frac{\sqrt{5}+1}{2} \left(\left\lceil\frac{q}{k}\right\rceil+1 \right). \]
To guarantee \eqref{1-12}, it is sufficient to require that
\begin{equation*}
	\delta_{3k} < \frac{-(\phi+1) + \sqrt{\phi^2+6\phi+1}}{2\phi} = \frac{2}{\sqrt{\phi^2 + 6\phi +1} + \phi+1}.
\end{equation*}
The right-hand side above is the positive root in $[0,1]$ of the quadratic equation $\phi \delta_{3k}^2 + \left( \phi + 1 \right)\delta_{3k} - 1 = 0$.
From the above analysis, we immediately obtain the following result.
\begin{corollary}
Let $x^* \in \mathbb{R}^n$ be the solution of the problem \eqref{problem} and $\eta := y-Ax^*$.
Let $q \geq 2 k$ be a positive integer number.
Suppose the restricted isometry constant of $A$ satisfy
\begin{equation} \label{1-16}
	\delta_{3k} < \frac{2}{\sqrt{\phi^2 + 6\phi +1} + \phi+1},
\end{equation}
where
\[ \phi = \frac{\sqrt{5}+1}{2}\left(\left\lceil\frac{q}{k}\right\rceil+1\right) . \]
Then the sequence $\left\{x^{p}\right\}$ generated by PGOT satisfies that
\begin{equation*}
	\left\|x^{p+1}-x^{*}\right\|_{2} \leq \rho^{p}\left\|x^{0}-x^{*}\right\|_{2}+\frac{\tau}{1-\rho}\|\eta\|_{2},
\end{equation*}
where $\rho<1$ and $\tau$ are given by \eqref{1-8} and \eqref{1-9}, respectively.
\end{corollary}

The bound \eqref{1-16} depends only on the given integer number $q$.
It is easy to verify that, for instance, $\delta_{3k} < 0.1517$ when $q = 2k$, $\delta_{3k} < 0.1211$ when $2k < q \le 3k$, and $\delta_{3k} < 0.1009$ when $3k < q \le 4k$.

Theorem 1 demonstrates how far the iterate point $x^{p+1}$ generated by PGOT is from the solution $x^*$ of the problem \eqref{problem}.
It shows that the bound of $\norm{x^{p+1} - x^*}_2$ depends on the value $\norm{\eta}_2 = \norm{y-Ax^*}_2$.
In many practical situations, for example in sparse signal recovery, $y$ are the linear measurements of the signal $x^*$.
In this case, $\norm{\eta}_2 = \norm{y-Ax^*}_2$ is the measurement error which is very small.
In particular, $\norm{\eta}_2=0$ when measurements are accurate.
In such practical problems, our error bounds established in Theorem 1 and Corollary 1 imply that $x^{p+1}$ generated by algorithms would approach to or even equal to $x^*$.
See the discussion below in more detail.

\subsubsection{Application to sparse signal recovery} \label{section3-1}
Let $x^*$ be a $k$-sparse signal to recover.
To recover $x^*$, we take the signal measurements $y := Ax^*+\eta$ with a measurement matrix $A$, where $\eta = y-Ax^*$ denotes the measurement error which is small.
Recovering $x^*$ from the measurements $y$ can be exactly modeled as the optimization problem \eqref{problem}.
From the results in Section \ref{section3}, we immediately obtain the next result concerning sparse signal recovery.

\begin{theorem} \label{PGOTcom}
Let $y:=A x^*+\eta$ be the measurements of the $k$-sparse signal $x^* \in \mathbb{R}^{n}$ with measurement error $\eta$.
Let $q \geq 2 k$ be a positive integer number.
Suppose the restricted isometry constant of measurement matrix A satisfies one the following conditions:
\begin{enumerate}[(i)]
\item $\delta_{3 k}<\alpha^*$, where $\alpha^* \in (0,1)$ is the unique real root of \eqref{alpha},
\item $\delta_{3 k}$ satisfies \eqref{1-16}.
\end{enumerate}
Then the sequence generated by PGOT satisfies
\begin{equation} \label{1-17}
\left\|x^{p+1}-x^*\right\|_{2} \leq \rho^p \left\|x^{0}-x^*\right\|_{2}+\frac{\tau}{1-\rho} \left\|\eta\right\|_{2},
\end{equation}
where $\rho$ and $\tau$ are the same as \eqref{1-8} and \eqref{1-9}, respectively.
In particular, if the measurements are accurate, i.e., $y=Ax^*$, then the sequence $\{x^p\}$generated by PGOT converges to $x^*$.
\end{theorem}

From \eqref{1-17}, we see that when the measurements are accurate enough, i.e., $\|\eta\|_2$ is sufficient small, then $x^{p+1} \approx x^*$.
This means the $x^{p+1}$ is a high-quality reconstruction of $x^*$.

\subsection{Main results for RPGOT}
Before analyzing the PGROT, we introduce the following lemma.
\begin{lemma} \emph{\cite{zhao2020optimal}} \label{lem5}
Let $x^* \in \mathbb{R}^n$ be the solution to the problem \eqref{problem} and $\eta = y-Ax^*$ be the error.
Denote $S = \supp(x^*)$ and $S^{p+1} = \supp(x^{p+1})$.
Let $u^p$ and $w^p$ be the vector defined as in PGROT, and $w^* \in \{0,1\}^n$ be a binary $k$-sparse vector such that $S \subseteq \supp(w^*)$.
Then
\begin{align*}
	\left\|\left(x^*-u^{p} \otimes w^{p}\right)_{S \cup S^{p+1}}\right\|_{2} \le & \sqrt{\frac{1+\delta_{k}}{1-\delta_{2 k}}} \left\|\left(x^*-u^{p}\right) \otimes w^*\right\|_{2} + \frac{2}{\sqrt{1-\delta_{2 k}}}\left\|\eta\right\|_{2} \\
	& + 2\sqrt{\frac{1+\delta_{k}}{1-\delta_{2 k}}} \left\|\mathcal{H}_{k}\left(u^{p}-x^*\right)\right\|_{2}  .
\end{align*}
\end{lemma}

\begin{theorem} \label{theorem3}
Let $x^* \in \mathbb{R}^n$ be the solution to the problem \eqref{problem} and $\eta = y-Ax^*$.
Let $q \geq 2 k$ be a positive integer number.
Suppose the restricted isometry constant of matrix $A$ satisfies
\begin{equation} \label{thm2-1}
	\delta_{3k} < \beta^*,
\end{equation}
where $\beta^*$ is the unique real root of the equation
\begin{equation} \label{beta}
	9\left(\left\lceil\frac{q}{k}\right\rceil+1\right)^{2} \beta^{3}+9 \left(\left\lceil\frac{q}{k}\right\rceil+1\right)^{2} \beta^{2}+\frac{2}{7+3 \sqrt{5}} \beta-\frac{2}{7+3 \sqrt{5}} = 0
\end{equation}
in $(0,1)$.
Then the sequence $\left\{x^{p}\right\}$ generated by PGROT satisfies
\begin{equation*}
\norm*{x^{p+1}-x^*}_2 \le \overline{\rho}^p \norm*{x^0-x^*}_2 + \frac{\overline{\tau}}{1-\overline{\rho}} \norm{\eta}_2,
\end{equation*}
where
\begin{equation} \label{thm2-2}
	\overline{\rho} = 3 \left(\frac{\sqrt{5}+1}{2}\right)^2 \left(\left\lceil\frac{q}{k}\right\rceil+1\right) \delta_{3k} \sqrt{\frac{1+\delta_{k}}{1-\delta_{2 k}}} < 1
\end{equation}
is ensured under \eqref{thm2-1}, and the constant $\overline{\tau}$ is given as
\begin{equation} \label{thm2-3}
	\overline{\tau} = \left(\frac{\sqrt{5}+1}{2}\right)^2 \frac{3\left(\left\lceil\frac{q}{k}\right\rceil+1\right) \left(1+\delta_k\right)}{\sqrt{1-\delta_{2 k}}} +\frac{\sqrt{5}+1}{\sqrt{1-\delta_{2 k}}} .
\end{equation}
In particular, if $\eta=0$, then the sequence $\left\{x^{p}\right\}$ generated by PGROT converges to $x^*$.
\end{theorem}
\emph{Proof.}
Let $x^{p+1}, u^p$ and $\overline{w}^p$ be defined in PGROT.
Denote by $S = \supp(x^*)$.
Note that $S^{p+1} = \supp(x^{p+1}) = \supp({\cal H}_k (\overline{w}^p \otimes u^p))$.
By Lemma \ref{hard}, we have
\begin{equation} \label{1-11}
\norm*{x^*-x^{p+1}}_2 = \left\| x^* - \mathcal{H}_{k}\left(\overline{w}^{p} \otimes u^{p}\right) \right\|_2 \le \frac{\sqrt{5}+1}{2} \norm*{(x^*-\overline{w}^{p} \otimes u^{p})_{S \cup S^{p+1}}}_2 .
\end{equation}
Note that $w^*$ is a $k$-sparse binary vector satisfying $\supp(x^*) \subseteq \supp(w^*)$.
By Lemma \ref{lem5}, we obtain
\begin{align}
& \left\|\left(x^*-u^{p} \otimes \overline{w}^{p}\right)_{S \cup S^{p+1}}\right\|_{2} \nonumber \\
& \le \sqrt{\frac{1+\delta_{k}}{1-\delta_{2 k}}} \left( \left\|\left(x^*-u^{p}\right) \otimes w^* \right\|_{2} + 2\left\|\mathcal{H}_{k}\left(u^{p}-x^*\right)\right\|_{2} \right) + \frac{2}{\sqrt{1-\delta_{2 k}}}\left\|\eta\right\|_{2} \nonumber \\
& \le 3 \sqrt{\frac{1+\delta_{k}}{1-\delta_{2 k}}} \norm*{x^*-u^p}_2+\frac{2}{\sqrt{1-\delta_{2 k}}}\left\|\eta\right\|_{2},
\end{align}
where the last inequality follows from
\[ \left\|\left(x^*-u^{p}\right) \otimes w^*\right\|_{2} \le \left\|\mathcal{H}_{k}\left(u^{p}-x^*\right)\right\|_{2} \le \norm*{x^*-u^p}_2. \]
Based on \eqref{1-2}, we have
\begin{equation} \label{2-2}
\norm*{x^*-u^p}_2 \le \frac{\sqrt{5}+1}{2}(t+1) \delta_{3 k}\left\|x^*-x^{p}\right\|_{2}+\frac{\sqrt{5}+1}{2}(t+1) \sqrt{1+\delta_{k}}\left\|\eta\right\|_{2} ,
\end{equation}
where $t = \left\lceil\frac{q}{k}\right\rceil$.
Combining \eqref{1-11} - \eqref{2-2} yields
\begin{align} \label{2-1}
\left\|x^*-x^{p+1}\right\|_{2} \le \overline{\rho} \left\|x^*-x^{p}\right\|_{2} + \overline{\tau} \norm*{\eta}_2
\end{align}
where $\overline{\rho}$ and $\overline{\tau}$ are given by \eqref{thm2-2} and \eqref{thm2-3}, respectively.
We now prove that \eqref{thm2-1} implies $\overline{\rho} < 1$.
Due to the fact $\delta_{k} \le \delta_{2k} \le \delta_{3k}$, to guarantee that $\overline{\rho} < 1$, it is sufficient to require
\begin{equation} \label{1-14}
3 \left(\frac{\sqrt{5}+1}{2}\right)^2 (t+1) \delta_{3k} \sqrt{\frac{1+\delta_{3k}}{1-\delta_{3k}}} < 1 ,
\end{equation}
which, by squaring both sides and rearranging terms, is equivalent to $g(\delta_{3k})<0$ where
\[ g\left(\delta_{3 k}\right) = 9 (t+1)^2 \delta_{3k}^3 + 9 (t+1)^2 \delta_{3k}^2 + \frac{2}{7+3\sqrt{5}} \delta_{3k} - \frac{2}{7+3\sqrt{5}}. \]
The gradient of $g(\delta_{3k})$ is given as
\[ \nabla g\left(\delta_{3 k}\right) = 27(t+1)^2 \delta_{3 k}^2 + 18 (t+1)^2 \delta_{3k} + \frac{2}{7+3\sqrt{5}}, \]
which is positive over the interval $[0,1]$.
This together with
\[ g(0) = - \frac{2}{7+3\sqrt{5}} < 0, ~ g(1) = 18 (t+1)^2 > 0 , \]
implies that there exists a unique real positive root $\beta^* \in (0,1)$ satisfying $g(\beta^*)=0$.
Therefore, the condition $\delta_{3k} < \beta^*$ guarantees the inequality \eqref{1-14}, and thus ensures that $\overline{\rho} < 1$.
Thus it follows from \eqref{2-1} that
\[ \norm*{x^{p+1}-x^*}_2 \le \overline{\rho}^p \norm*{x^0-x^*}_2 + \frac{\overline{\tau}}{1-\overline{\rho}} \norm{\eta}_2. \]

When $\eta=0$, the relation above implies that $\left\|x^*-x^{p+1}\right\|_{2} \leq \overline{\rho}^{p}\left\|x^*-x^{0}\right\|_{2} \to 0$ as $p \to \infty$.
Therefore, the sequence $\{x^p\}$ generated by RPGOT in this case converges to the solution $x^*$ of \eqref{problem}.
\hfill $\Box$

Similar to the discussion in the end of Section 3.1, an explicit bound of $\delta_{3k}$ for PGROT can be given.
Since $\sqrt{\frac{1+\delta_{3 k}}{1-\delta_{3 k}}} < \frac{1+\delta_{3 k}}{1-\delta_{3 k}}$, a sufficient condition for \eqref{1-14} is
\[ 3 \left(\frac{\sqrt{5}+1}{2}\right)^2 \left(\left\lceil\frac{q}{k}\right\rceil+1 \right) \delta_{3k} \frac{1+\delta_{3k}}{1-\delta_{3k}} < 1, \]
which is equivalent to
$$\psi \delta_{3 k}^{2}+(\psi+1) \delta_{3 k}-1 < 0, $$
where
\[ \psi = 3\left(\frac{\sqrt{5}+1}{2}\right)^{2}\left(\left\lceil\frac{q}{k}\right\rceil+1\right). \]
By the same analysis in Section 3.1, we immediately have the next corollary.

\begin{corollary} \label{corollary2}
Let $x^* \in \mathbb{R}^n$ be the solution to the problem \eqref{problem} and $\eta := y-Ax^*$.
Let $q \geq 2 k$ be a positive integer number.
Suppose the restricted isometry constant of matrix $A$ satisfies
\[ \delta_{3k} < \frac{2}{\sqrt{\psi^2 + 6\psi +1} + \psi+1}, \]
where
\[ \psi = 3\left(\frac{\sqrt{5}+1}{2}\right)^{2}\left(\left\lceil\frac{q}{k}\right\rceil+1\right) . \]
Then the sequence $\left\{x^{p}\right\}$ generated by PGROT satisfies
\begin{equation*}
\norm*{x^{p+1}-x^*}_2 \le (\overline{\rho})^p \norm*{x^{0}-x^*}_2 + \frac{\overline{\tau}}{1-\overline{\rho}} \norm{\eta}_2,
\end{equation*}
where $\overline{\rho}$ and $\overline{\tau}$ are given as \eqref{thm2-2} and \eqref{thm2-3}, respectively.
\end{corollary}

Similar to Corollary 1, we may apply the above result (Theorem \ref{theorem3} and Corollary \ref{corollary2}) to the scenario of sparse signal recovery via compressed sensing for which $\eta = y-Ax^*$ is very small, and thus $x^p \approx x^*$ when $p$ is large enough.
That is, the $x^p$ generated by PGROT is a quality approximation to the signal.

\subsection{Main result for RPGOTP}
The error bound for the solution of \eqref{problem} via PGROTP algorithm can be also established.
The next lemma concerning a property of pursuit step is useful in this analysis.
\begin{lemma} \emph{\cite{zhao2020optimal}} \label{lem8}
Let $x^* \in \mathbb{R}^n$ be the solution to the problem \eqref{problem} and $\eta = y-Ax^*$.
The vector $u \in \mathbb{R}^n$ is an arbitrary $k$-sparse vector.
Then the optimal solution of the pursuit step
\[ z^{*}=\arg \min _{z}\left\{\|y-A z\|_{2}^{2}: ~ \operatorname{supp}(z) \subseteq \operatorname{supp}(u)\right\} \]
satisfies that
\[ \left\|z^{*}-x^*\right\|_{2} \leq \frac{1}{\sqrt{1-\left(\delta_{2 k}\right)^{2}}}\|x^*-u\|_{2}+\frac{\sqrt{1+\delta_{k}}}{1-\delta_{2 k}}\|\eta\|_{2} . \]
\end{lemma}
The main result for PGROTP algorithm is given as follows.
\begin{theorem}
Let $x^* \in \mathbb{R}^n$ be the solution to the problem \eqref{problem} and $\eta = y-Ax^*$.
Let $q \geq 2 k$ be a positive integer number.
Suppose the restricted isometry constant of matrix $A$ satisfies
\begin{equation} \label{gamma}
	\delta_{3k} < \frac{1}{3\left(\frac{\sqrt{5}+1}{2}\right)^{2}\left(\left\lceil\frac{q}{k}\right\rceil+1\right) + 1}.
\end{equation}
Then the sequence $\left\{x^{p}\right\}$ generated by PGROTP satisfies
\begin{equation} \label{thm3-3}
	\norm*{x^{p+1}-x^*}_2 \le (\rho^*)^p \norm*{x^0-x^*}_2 + \frac{\tau^*}{1-\rho^*} \norm{\eta}_2,
\end{equation}
where
\begin{equation} \label{thm3-1}
\rho^* = \left(\frac{\sqrt{5}+1}{2}\right)^{2} \frac{3\left(\left\lceil\frac{q}{k}\right\rceil+1\right)\delta_{3 k}}{1-\delta_{3k}} < 1
\end{equation}
is guaranteed by \eqref{gamma}, and the constant $\tau^*$ is given as
\begin{equation} \label{thm3-2}
\tau^* =  \frac{3\left( \frac{\sqrt{5}+1}{2} \right)^2\left(\left\lceil\frac{q}{k}\right\rceil+1\right)(1+\delta_k)+\sqrt{5}+1}{(1-\delta_{2k}) \sqrt{1+\delta_{2k}}} + \frac{\sqrt{1+\delta_k}}{1-\delta_{2k}}.
\end{equation}
In particular, when $\eta=0$, the sequence $\left\{x^{p}\right\}$ generated by PGROTP converges to $x^*$.
\end{theorem}
\emph{Proof.}
The PGROTP can be regarded as a combination of PGROT with a pursuit step.
Denote $\overline{x}^{p+1}$ as the intermediate vector generated by PGROT.
Based on the analysis of PGROT algorithm, we have
\[ \norm{\overline{x}^{p+1} -x^*}_2 \le \overline{\rho} \left\|x^*-x^{p}\right\|_{2} + \overline{\tau} \norm*{\eta}_2, \]
where $\overline{\rho}$ and $\overline{\tau}$ are the same as \eqref{thm2-2} and \eqref{thm2-3}, respectively.
By Lemma \ref{lem8}, we immediately have that
\[ \left\|x^*-x^{p+1}\right\|_{2} \le \frac{1}{\sqrt{1-\left(\delta_{2 k}\right)^{2}}} \norm*{\overline{x}^{p+1} -x^*}_2 + \frac{\sqrt{1+\delta_{k}}}{1-\delta_{2 k}} \norm*{\eta}_2 . \]
As $\delta_k \le \delta_{2k} \le \delta_{3k}$, combining two inequalities above yields
\[ \left\|x^*-x^{p+1}\right\|_{2} \le \rho^* \left\|x^*-x^{p}\right\|_{2} + \tau^* \norm*{\eta}_2, \]
where $\rho^*$ and $\tau^*$ are given in \eqref{thm3-1} and \eqref{thm3-2}, respectively.
To guarantee $\rho^* < 1$, i.e.,
\[ \left(\frac{\sqrt{5}+1}{2}\right)^{2} \frac{3\left(\left\lceil\frac{q}{k}\right\rceil+1\right)\delta_{3 k}}{1-\delta_{3k}} < 1, \]
it is sufficient to require that
\[ \delta_{3k} < \frac{1}{3\left(\frac{\sqrt{5}+1}{2}\right)^{2} \left(\left\lceil\frac{q}{k}\right\rceil+1\right) + 1}, \]
which is exactly the assumption \eqref{thm3-3} of the theorem.

If $\eta=0$, the sequence $\left\{x^{p}\right\}$ generated by PGROTP converges to $x^*,$   since in this case, \eqref{thm3-3} is reduced to $\norm*{x^{p+1} - x^*}_2 \le (\rho^*)^p \norm*{x^0-x^*}_2\to 0$ as $p \rightarrow \infty. $
\hfill $\Box$

For sparse signal recovery, similar comments to that of Section \ref{section3-1} can be made to PGROTP.
The discussion is omitted here.
Before we close this section, we list a few RIC conditions in terms of $\delta_{3k}$ for the proposed algorithms with different $q$, i.e., $q=2k,3k,4k$.
The results shown in Table 1 are derived based on \eqref{alpha}, \eqref{beta} and \eqref{gamma} for the $q$ as given above, respectively.
It is worth mentioning that, when $q=n$, the partial gradient ${\cal H}_q(\nabla f(x))$ becomes the full gradient $\nabla f(x)$, and the algorithms in this paper are reduced to the optimal $k$-thresholding (OT), relaxed optimal $k$-thresholding (ROT) and relaxed optimal $k$-thresholding pursuit (ROTP) algorithm, respectively.
The sufficient conditions for the convergence of these algorithms were studied in \cite{zhao2020optimal,zhao2020analysis}.

\begin{table}[H]\centering
\caption{The upper bounds of $\delta_{3k}$ for several different $q$}
\begin{tabular}{|c|c|c|c|}
\hline
The value of $q$ & PGOT                   & RPGOT                  & RPGOTP                 \\ \hline
$q=2k$           & $\delta_{3k} < 0.1729$ & $\delta_{3k} < 0.0407$ & $\delta_{3k} < 0.0407$ \\ \hline
$2k < q \le 3k$           & $\delta_{3k} < 0.1348$ & $\delta_{3k} < 0.0308$ & $\delta_{3k} < 0.0308$ \\ \hline
$3k < q \le 4k$ & $\delta_{3k} < 0.1106$ & $\delta_{3k} < 0.0248$ & $\delta_{3k} < 0.0248$ \\ \hline
\end{tabular}
\captionsetup{justification=centering}
\label{table1}
\end{table}

\section{Numerical Experiments}
Simulations via synthetic data are carried out to demonstrate the numerical performance of the PGROTP which is the main implementable algorithm proposed in this paper.
We test the algorithm from three aspects: objective reduction, average number of iterations required for solving the problem \eqref{problem}, and success frequency in vector reconstruction.
The PGROTP with $q=k, 2k, 3k$ and $n$ are tested and compared.
The measurement matrices used in experiments are Gaussian random matrices whose entries follow standard normal distribution ${\cal N}(0,1)$.
For sparse vectors, their entries also follow the ${\cal N}(0,1)$ and the position of nonzero entries of the vector follows the uniform distribution.
All involved convex optimization problems were solved by CVX developed by Grant and Boyd \cite{grant2017cvx} with solver 'Mosek' \cite{andersen2000mosek}.

\subsection{Objective reduction}
\begin{figure}[!htp]
\centering
\subfigure[Sparsity level $k=162$]{
\includegraphics[width=0.45\textwidth]{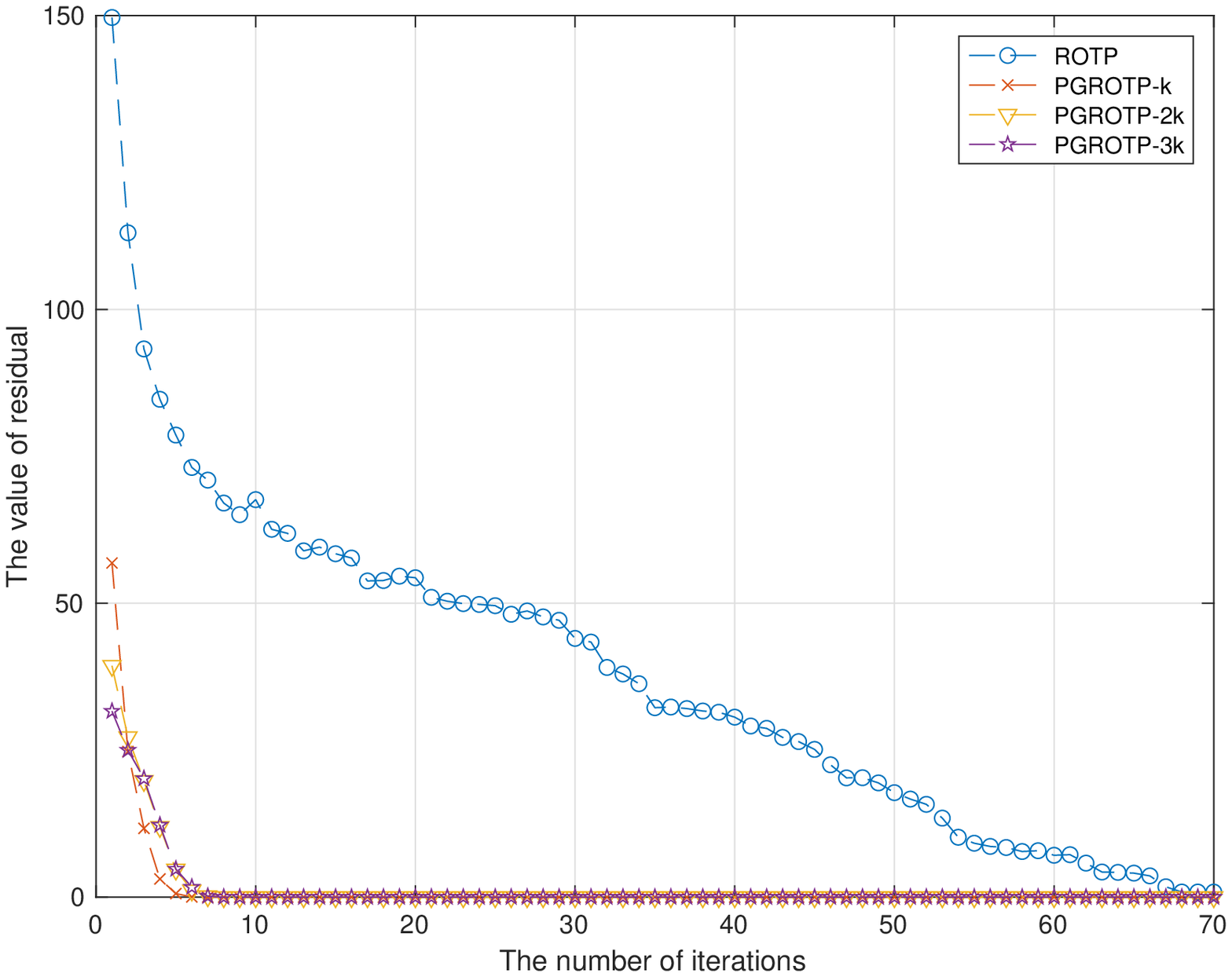}}
\subfigure[Sparsity level $k=197$]{
\includegraphics[width=0.45\textwidth]{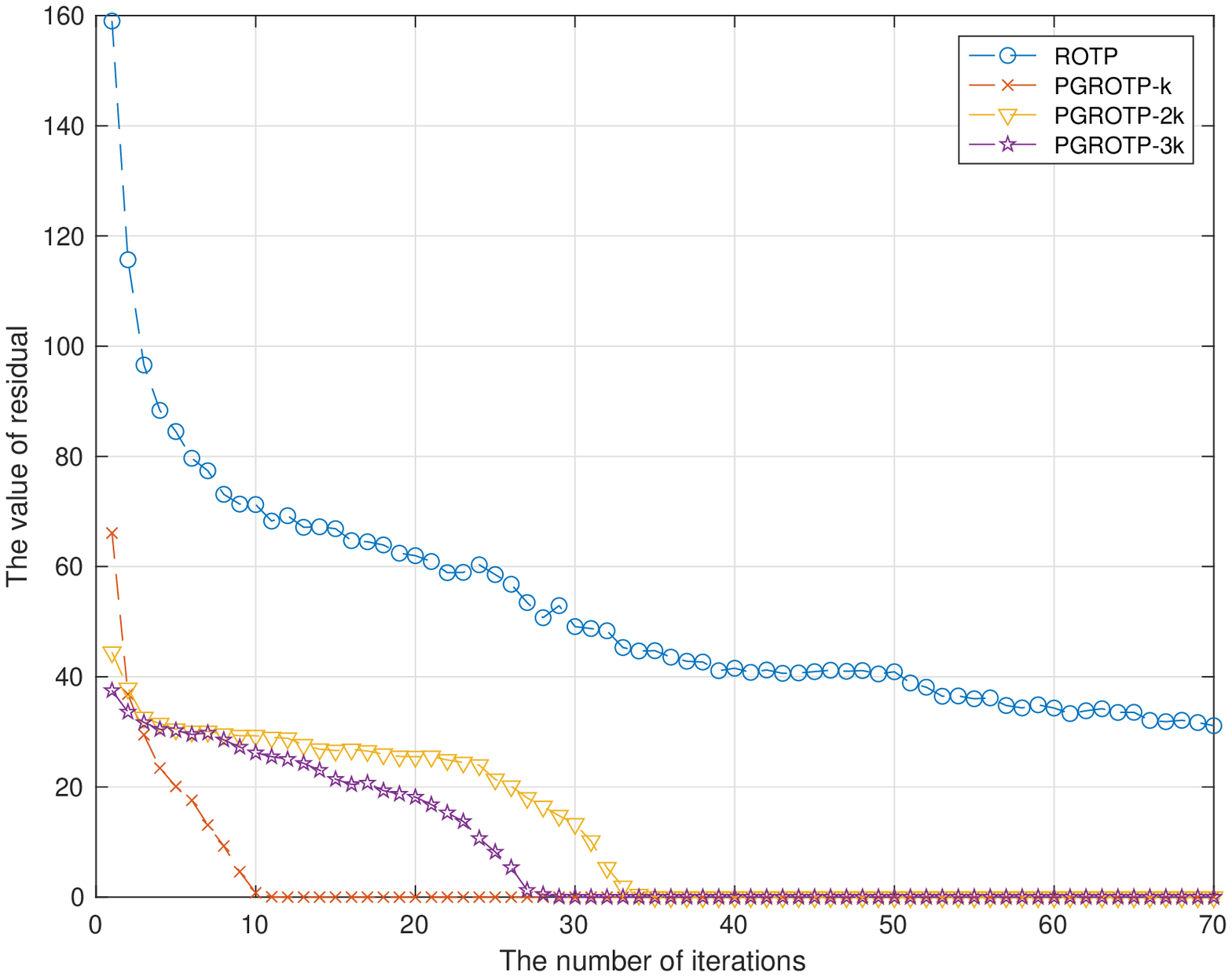}}
\caption{Objective change in the course of iterations for PGROTP with different $q$, i.e., $q = k, 2k, 3k, n$.}
\label{graph1}
\end{figure}

This experiment is used to investigate the objective-reduction performance of the PGROTP with different $q$, including $q= k,2k,3k$ and $n$.
We set $A \in \mathbb{R}^{500 \times 1000}$ and $y=Ax^*$, where $x^*$ is a generated $k$-sparse vector.
Thus $x^*$ is a global solution of the problem \eqref{problem}.
Fig. \ref{graph1} records the changes of the objective value $\norm{y-Ax}_2$ in the course of algorithm up to 70 iterations.
Fig. \ref{graph1} (a) and Fig. \ref{graph1} (b) include the results for the sparsity level $\norm{x^*}_0=162$ and $197$, respectively.
It can be seen that PGROTP is able to reduce the objective value during iterations.
Moreover, this experiment also indicates that the optimal $k$-thresholding methods with partial gradients often perform better in objective reduction than the ones using full gradients.

\subsection{Number of iterations}
Experiments were also performed to demonstrate the average number of iterations needed for PGROTP to solve the sparse optimization problems from the sparse vector reconstruction.
The vector dimension is fixed to be 1000, and the size of the measurement matrix is $m \times 1000$, where $m$ takes the following a few different values: $m = 300, 400, 500, 600$.
The measurements $y = Ax^*$ are accurate, where $x^*$ is the sparse vector to recover.
In this experiment, if the iterate $x^p$ satisfies the criterion
\begin{equation} \label{stopcriteria}
{\norm{x-x^*}_2}/{\norm{x^*}_2} \le 10^{-3},
\end{equation}
then the algorithm terminates and the number of iterations $p$ is recorded.
If the algorithm within 50 iterations cannot meet the criterion \eqref{stopcriteria}, then the algorithm stops, and the number of iterations performed is recorded as 50.
For each given sparsity level, the average number of iterations is obtained by attempting 50 trials.
\begin{figure}
\centering
\subfigure[$m=300$]{
\includegraphics[width=0.45\textwidth]{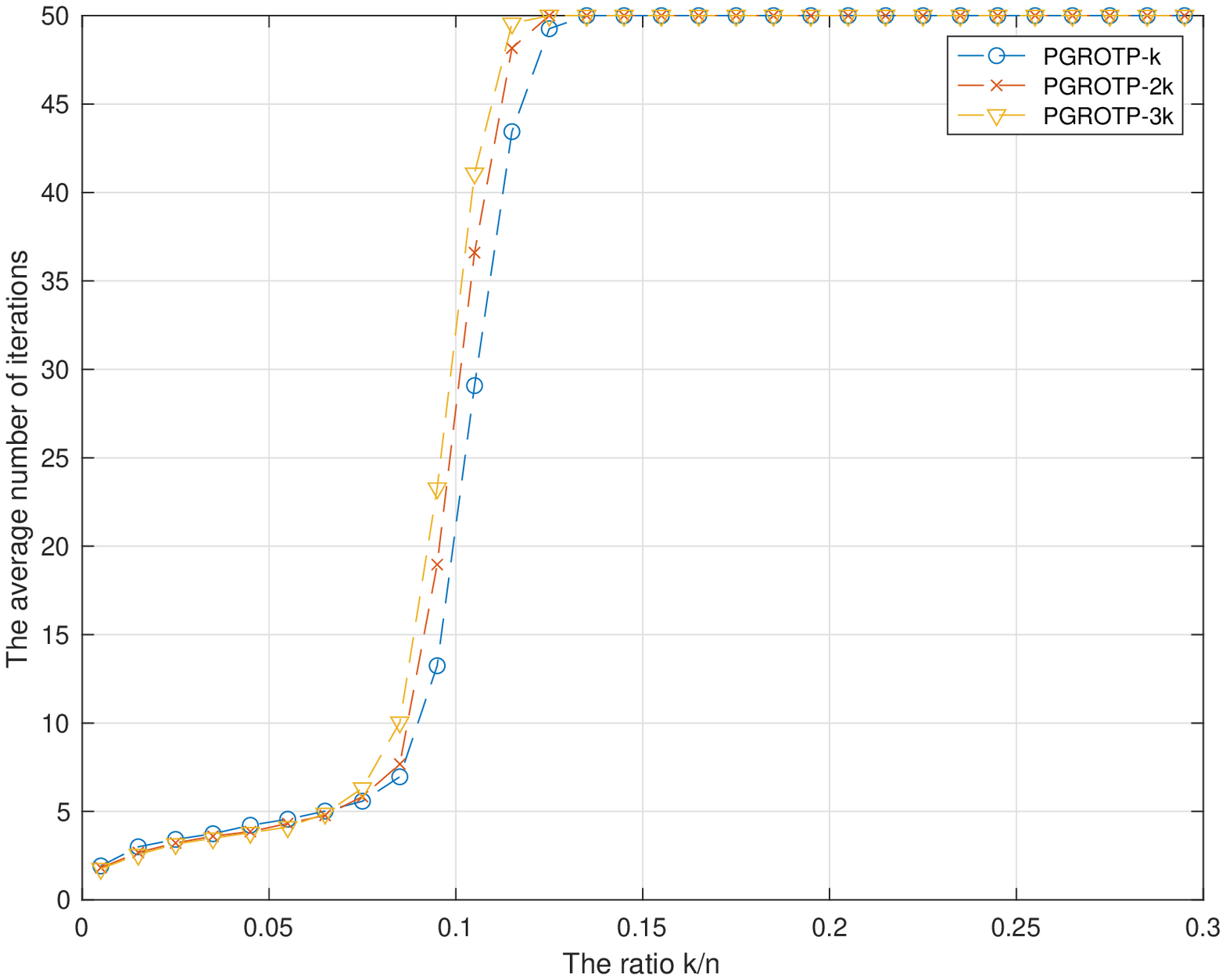}}
\subfigure[$m=400$]{
\includegraphics[width=0.45\textwidth]{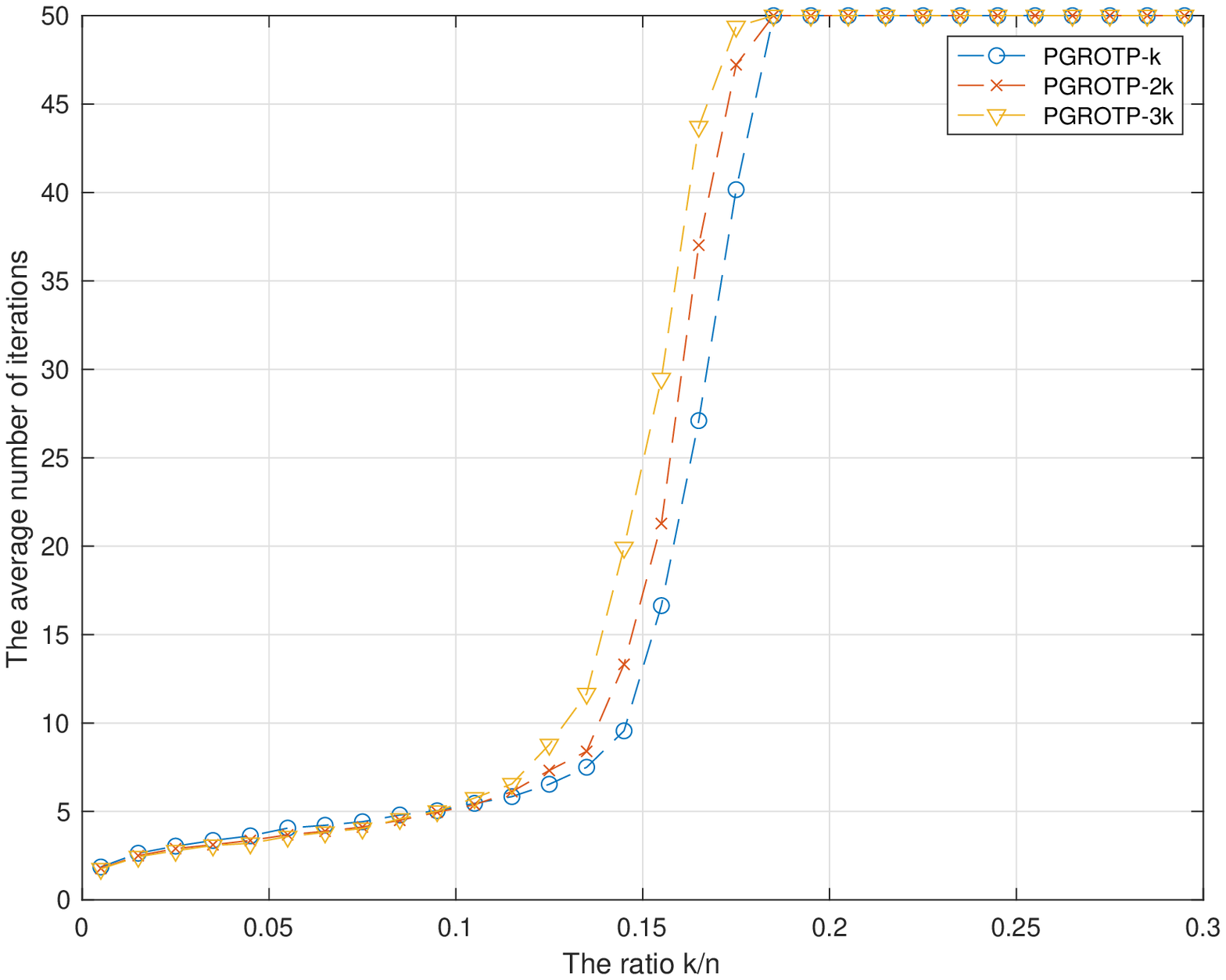}}
\subfigure[$m=500$]{
\includegraphics[width=0.45\textwidth]{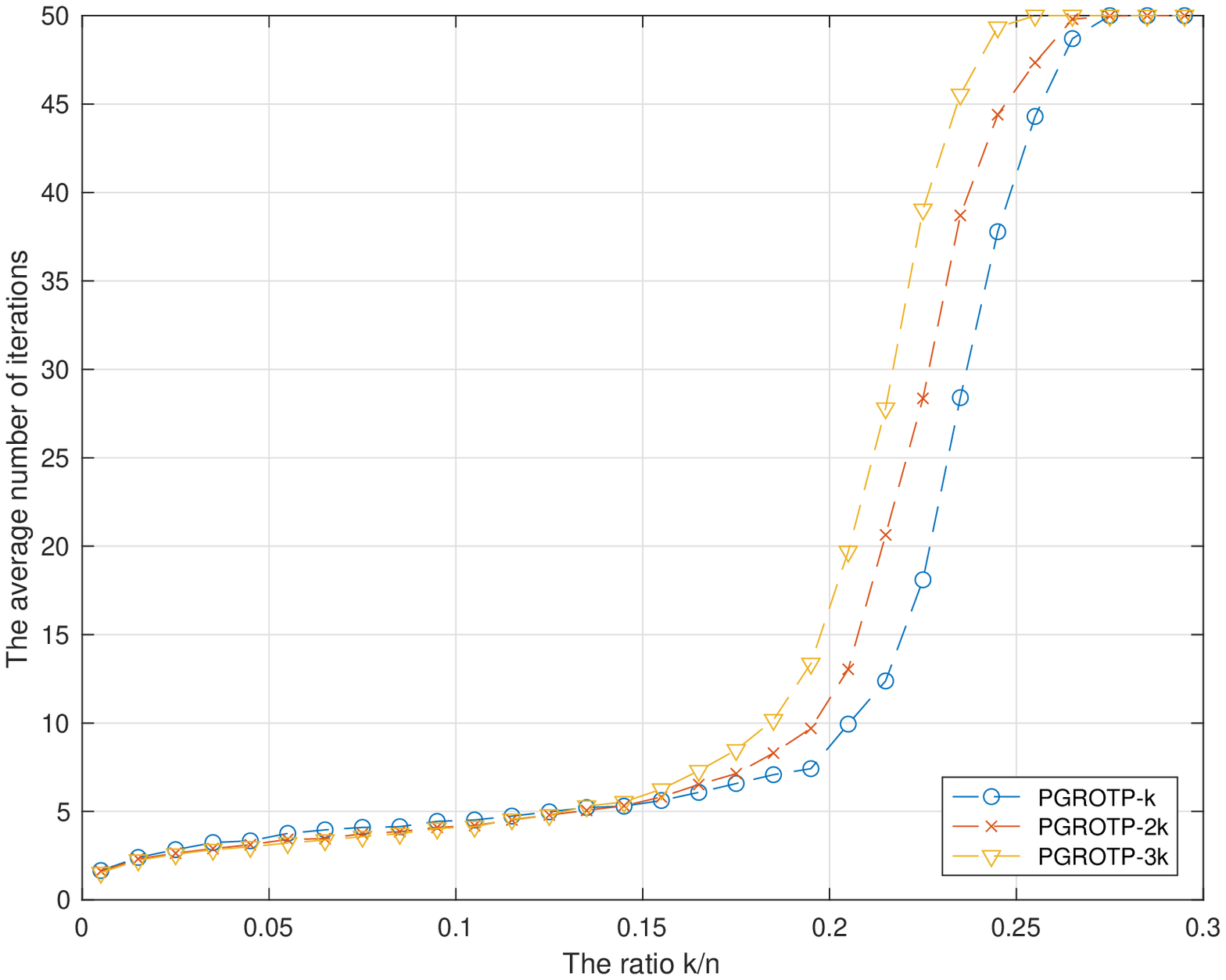}}
\subfigure[$m=600$]{
\includegraphics[width=0.45\textwidth]{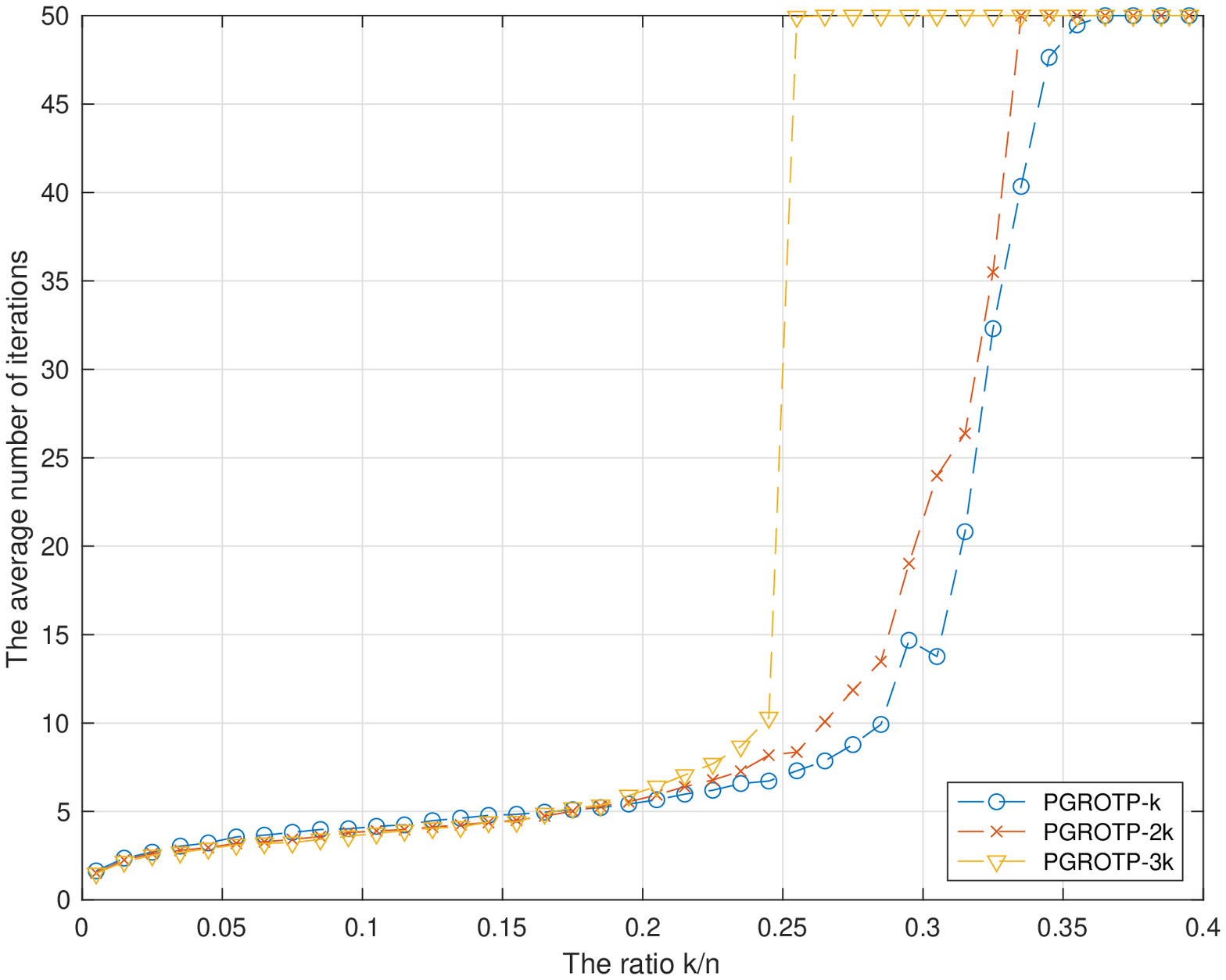}}
\caption{Comparison of the average number of iterations required by PGROTP with different $q$}
\label{graph2}
\end{figure}
The outcomes are shown in Fig. \ref{graph2} which indicate that the required iterations of PGROTP for vector reconstructions are usually low when the sparsity level of $x^*$ is low, and that the number of iterations which is required for solving the problem increases as the sparsity level increases.
This figure also shows that which is for a given sparsity level, the more measurements are required, the lower the average number of iterations are needed by the PGROTP to meet the reconstruction criterion \eqref{stopcriteria}.

\subsection{Sparse signal recovery}

\begin{figure}
\centering
\subfigure[Exact measurements]{
\includegraphics[width=0.45\textwidth]{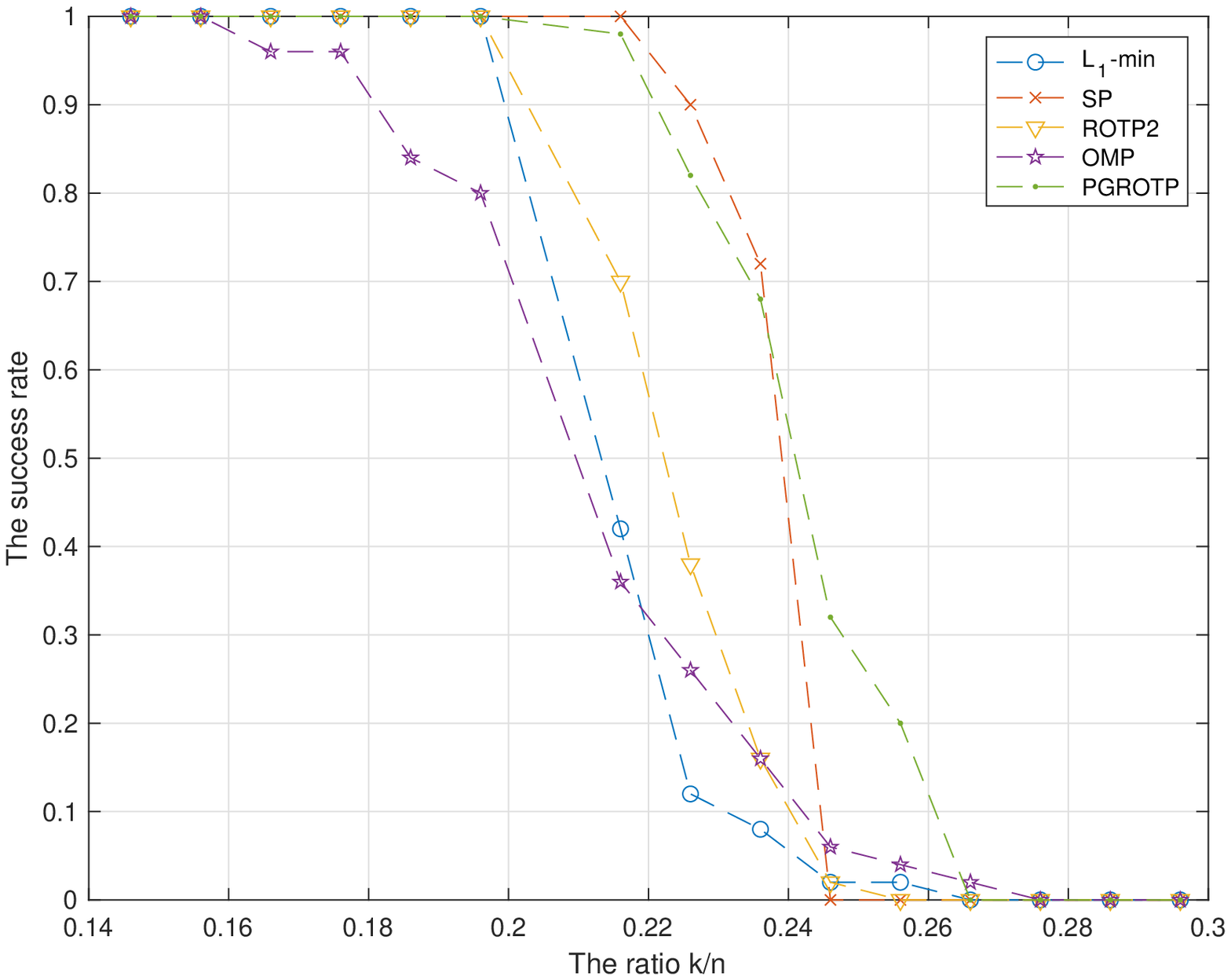}}
\subfigure[Inexact measurements]{
\includegraphics[width=0.45\textwidth]{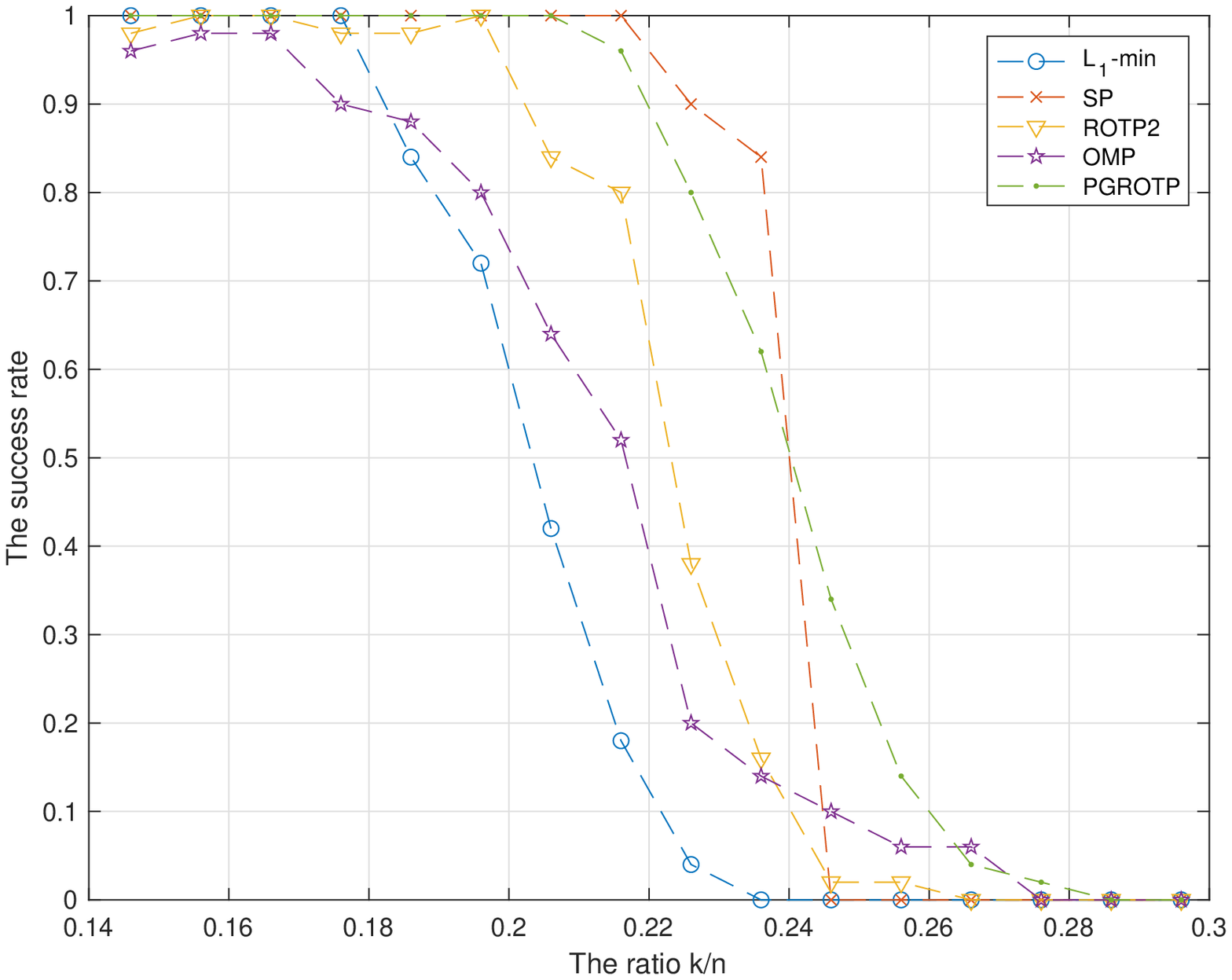}}
\caption{Comparisons of success rates of sparse signal reconstruction between several algorithms via Gaussian random matrices}
\label{graph3}
\end{figure}

Simulations were also carried out to compare the success rates of the PGROTP algorithm in sparse vector reconstruction with several existing algorithms, such as $\ell_1$-minimization, subspace pursuit (SP), orthogonal matching pursuit (OMP) and ROTP2 (in \cite{zhao2020optimal,zhao2020analysis}).
The size of the measurement matrix is still ${500 \times 1000}$.
For every given ratio of the sparsity level $k$ and $n$, the success rate of the algorithm is obtained by 50 random attempts.
In this experiment, SP, ROTP2 and PGROTP perform a total of 50 iterations, whereas OMP is performed $k$ iterations.
After performing the required number of iterations, the algorithm is counted as success if the condition \eqref{stopcriteria} is satisfied.
The success rates for accurate and inaccurate measurements are summarized in  Fig. \ref{graph3} (a) and (b), respectively.
The inaccurate measurements are given as $y=Ax^*+0.001\eta$, where $\eta$ is a standard Gaussian vector.
Compared with several existing algorithms, it can be seen that the PGROTP is robust and efficient for the sparse vector reconstruction in both noise and noiseless environment.

\section{Conclusions}
Motivated by the recent optimal $k$-thresholding technique, we proposed the partial-gradient-based optimal $k$-thresholding methods for solving a class of sparse optimization problems.
Under the restricted isometry property, we established a global error bound for the iterates produced by our algorithms.
Reduced to sparse signal recovery, our results claim that the proposed algorithms with $q$ satisfying $2k \le q < n$ are guaranteed to recover the sparse vector.
Numerical experiments demonstrate that the PGROTP algorithm is efficient for sparse vector reconstruction.
Although we focus on solving the specific model \eqref{problem} in this paper, it is not difficult to extend the framework of the proposed algorithms to the general model \eqref{prob}.
We leave this as a future work.


\begin{thebibliography}{999} \label{ref}

\bibitem{andersen2000mosek} Andersen, E. D., Andersen, K. D.: The MOSEK interior point optimizer for linear programming: an implementation of the homogeneous algorithm. High performance optimization. Springer, Boston, MA, 33, 197-232 (2000)

\bibitem{blumensath2008iterative} Blumensath, T., Davies, M.: Iterative hard thresholding for sparse approximation. J. Fourier Anal. Appl. 14, 629-654 (2008)

\bibitem{blumensath2009iterative} Blumensath, T., Davies, M.: Iterative hard thresholding for compressed sensing. Appl. Comput. Harmon. Anal. 27, 265-274 (2009)

\bibitem{blumensath2010normalized} Blumensath, T., Davies, M.: Normalized iterative hard thresholding: Guaranteed stability and performance. IEEE J. Sel. Top. Signal Process. 4, 298-309 (2010)

\bibitem{blumensath2012accelerated} Blumensath, T.: Accelerated iterative hard thresholding. Signal Process. 92, 752-756 (2012)

\bibitem{bouchot2016hard} Bouchot, J., Foucart, S., Hitczenki, P.: Hard thresholding pursuit algorithms: Number of iterations. Appl. Comput. Harmon. Anal. 41, 412-435 (2016)

\bibitem{boche2019compressed} Boche, H., Calderbank, R., Kutyniok, G., Vybiral, J.: Compressed Sensing and Its Applications. Springer, New York (2019)

\bibitem{cai2011orthognal} Cai, T., Wang, L.: Orthogonal matching pursuit for sparse signal recovery with noise. IEEE Transactions on Information theory. 57(7), 4680-4688 (2011)

\bibitem{candes2005decoding} Cand{\`e}s, E., Tao, T.: Decoding by linear programming. IEEE Trans. Inform. Theory. 51(12), 4203-4215 (2005)

\bibitem{candes2006robust} Cand{\`e}s, E., Romberg, J., Tao, T.: Robust uncertainty principles: Exact signal reconstruction from highly incomplete frequency information. IEEE Trans. Inform. Theory. 52(2), 489-509 (2006)

\bibitem{candes2008enhancing} Cand{\`e}s, E., Wakin, M., Boyd, S.: Enhancing sparsity by reweighted $\ell_1$-minimization. J. Fourier Anal. Appl. 14, 877-905 (2008)

\bibitem{chen2001atomic} Chen, S., Donoho, D., Saunders, M.: Atomic decomposition by basis pursuit. SIAM Review. 43(1), 129-159 (2001)

\bibitem{choi2017compressed} Choi, J., Shim, B., Ding, Y., Rao, B., Kim, D.: Compressed sensing for wireless communi- cations: Useful tips and tricks. IEEE Commun. Surveys $\&$ Tutorials. 19(3), 1527-1549, (2017)

\bibitem{dai2008subspace} Dai, W., Milenkovic, O.: Subspace pursuit for compressive sensing: Closing the gap between performance and complexity. ILLINOIS UNIV AT URBANA-CHAMAPAIGN (2008)

\bibitem{dai2009subspace} Dai, W., Milenkovic, O.: Subspace pursuit for compressive sensing signal reconstruction. IEEE Trans. Inform. Theory. 55, 2230-2249 (2009)

\bibitem{donoho1995denoising} Donoho, D.: De-noising by soft-thresholding. IEEE Trans. Inform. Theory. 41, 613-627 (1995)

\bibitem{eldar2012compressed} Eldar, Y., Kutyniok, G.: Compressed Sensing: Theory and Applications. Cambridge University Press, Cambridge (2012)

\bibitem{foucart2011hard} Foucart, S.: Hard thresholding pursuit: an algorithm for compressive sensing. SIAM J. Numer. Anal. 49(6), 2543-2563 (2011)

\bibitem{foucart2013mathematical} Foucart, S., Rauhut, H.: A Mathematical Introduction to Compressive Sensing. Birkh\"{a}user Basel (2013)

\bibitem{fornasier2008iterative} Fornasier, M., Rauhut, H.: Iterative thresholding algorithms. Appl. Comput. Harmon. Anal. 25(2), 187-208 (2008)

\bibitem{garg2009gradient} Garg, Rahul, and Rohit Khandekar. "Gradient descent with sparsification: an iterative algorithm for sparse recovery with restricted isometry property." Proceedings of the 26th annual international conference on machine learning. 2009.

\bibitem{huang2013soft} Huang, G., Wang, L.: Soft-thresholding orthogonal matching pursuit for efficient signal reconstruction. 2013 IEEE International Conference on Acoustics, Speech and Signal Processing. 2543-2547 (2013)

\bibitem{khanna2018iht} Khanna R, Kyrillidis A. IHT dies hard: Provable accelerated iterative hard thresholding. International Conference on Artificial Intelligence and Statistics. 188-198 (2018)

\bibitem{liu2016projected} Liu Y, Zhan Z, Cai J F, et al. Projected iterative soft-thresholding algorithm for tight frames in compressed sensing magnetic resonance imaging. IEEE transactions on medical imaging. 35(9). 2130-2140 (2016)

\bibitem{grant2017cvx} Grant, M., Boyd, S.: \textit{CVX: Matlab software for Disciplined Convex Programming}. Version 1.21, April 2017.

\bibitem{nan2020newton} Meng, N., Zhao, Y.-B.: Newton-step-based hard thresholding algorithms for sparse signal recovery. IEEE Trans. Signal Process. 68, 6594-6606 (2020)

\bibitem{nan2021newton} Meng N., Zhao Y.-B.: Newton-type optimal thresholding algorithms for sparse optimization problems. arXiv preprint arXiv:2104.02371, 2021.

\bibitem{needell2009cosamp} Needell, D., Tropp, J.: CoSaMP: Iterative signal recovery from incomplete and inaccurate samples. Appl. Comput. Harmon. Anal. 26, 301-321 (2009)

\bibitem{needell2010signal} Needell, D., Vershynin, R.: Signal recovery from incomplete and inaccurate measurements via regularized orthogonal matching pursuit. IEEE J. Sel. Top. Signal Process. 4(2), 310-316 (2010)

\bibitem{patel2011sparse} Patel, V., Chellappa, R.: Sparse representations, compressive sensing and dictionaries for pattern recognition. The First Asian Conference on Pattern Recognition, IEEE. 325-329 (2011)

\bibitem{tropp2007signal} Tropp, J., Gilbert, A.: Signal recovery from random measurements via orthogonal mathcing pursuit. IEEE Trans. Inform. Theory. 53(12), 4655-4666 (2007)

\bibitem{yuan2014newton} Yuan, X.-T., Liu, Q.: Newton greedy pursuit: A quadratic approximation method for sparsity-constrained optimization. Proceedings of the IEEE Conference on Computer Vision and Pattern Recognition. 4122-4129 (2014)

\bibitem{zhao2018sparse}  Zhao, Y.-B.: Sparse Optimization Theory and Methods. CRC Press, Boca Raton, FL, (2018)

\bibitem{zhao2020optimal} Zhao, Y.-B.: Optimal $k$-thresholding algorithms for sparse optimization problems. SIAM J. Optim. 30(1), 31-55 (2020)

\bibitem{zhao2015new} Zhao, Y.-B., Ko$\check{\text{c}}$vara, M.: A new computational method for the sparsest solutions to systems of linear equations. SIAM J. Optim. 25(2), 1110-1134 (2015)

\bibitem{zhao2012reweighted} Zhao, Y.-B., Li, D.: Reweighted $\ell_1$-minimization for sparse solutions to underdetermined linear systems. SIAM J. Optim. 22(3), 1065-1088 (2012)

\bibitem{zhao2017constructing} Zhao, Y.-B., Luo, Z.-Q.: Constructing new reweighted $\ell_1$-algorithms for sparsest points of polyhedral sets. Math. Oper. Res. 42, 57-76 (2017)

\bibitem{zhao2020analysis} Zhao, Y.-B., Luo, Z.-Q.: Analysis of optimal thresholding algorithms for compressed sensing. Signal Process. 187, p. 108148 (2021)

\bibitem{zhao2020improved} Zhao, Y.-B., Luo, Z.-Q.: Improved RIP-Based Bounds for Guaranteed Performance of two Compressed Sensing Algorithms. arXiv:2007.01451v3, 2020.

\bibitem{zhou2021global} Zhou, S., Xiu, N., Qi, H.: Global and quadratic convergence of Newton hard-thresholding pursuit. arXiv:1901.02763v1 (2019)

\bibitem{zhou2020subspace} Zhou, S., Pan, L., Xiu, N.: Subspace Newton method for the $\ell_0 $-regularized optimization. arXiv:2004.05132, (2020)

\end{thebibliography}
\end{document}